\documentclass[12pt,leqno]{article}
\usepackage{amssymb}
\usepackage[mathscr]{eucal}
\usepackage{amsmath,amssymb,latexsym,theorem,bbm}
\usepackage{color}
\usepackage{appendix}

\newcommand{\SC}{\scriptstyle}

\newcommand{\CC}{\mathsf{C}}
\newcommand{\DD}{\mathsf{D}}
\newcommand{\NN}{\mathbb{N}}
\newcommand{\RR}{\mathbb{R}}
\newcommand{\ZZ}{\mathbb{Z}}

\newcommand{\bA}{{\boldsymbol{A}}}
\newcommand{\bC}{{\boldsymbol{C}}}
\newcommand{\bBB}{{\boldsymbol{B}}}
\newcommand{\bI}{{\boldsymbol{I}}}
\newcommand{\bm}{{\boldsymbol{m}}}
\newcommand{\bM}{{\boldsymbol{M}}}
\newcommand{\cM}{{\mathcal M}}
\newcommand{\bcM}{\boldsymbol{\cM}}
\newcommand{\bq}{{\boldsymbol{q}}}
\newcommand{\bu}{{\boldsymbol{u}}}
\newcommand{\bU}{{\boldsymbol{U}}}
\newcommand{\bv}{{\boldsymbol{v}}}
\newcommand{\bV}{{\boldsymbol{V}}}

\newcommand{\bx}{{\boldsymbol{x}}}
\newcommand{\bX}{{\boldsymbol{X}}}
\newcommand{\by}{{\boldsymbol{y}}}
\newcommand{\bW}{{\boldsymbol{W}}}
\newcommand{\balpha}{{\boldsymbol{\alpha}}}
\newcommand{\bgamma}{{\boldsymbol{\gamma}}}
\newcommand{\bxi}{{\boldsymbol{\xi}}}
\newcommand{\bmu}{{\boldsymbol{\mu}}}
\newcommand{\bPi}{{\boldsymbol{\Pi}}}
\newcommand{\bvare}{{\boldsymbol{\vare}}}
\newcommand{\bzero}{{\boldsymbol{0}}}

\newcommand{\cA}{{\mathcal A}}
\newcommand{\cB}{{\mathcal B}}
\newcommand{\cD}{{\mathcal D}}
\newcommand{\cF}{{\mathcal F}}
\newcommand{\cL}{{\mathcal L}}
\newcommand{\cP}{{\mathcal P}}
\newcommand{\cQ}{{\mathcal Q}}
\newcommand{\bcQ}{{\boldsymbol{\cQ}}}
\newcommand{\cR}{{\mathcal R}}
\newcommand{\bcR}{{\boldsymbol{\cR}}}
\newcommand{\cU}{{\mathcal U}}
\newcommand{\bcU}{\boldsymbol{\cU}}
\newcommand{\cX}{{\mathcal X}}
\newcommand{\bcX}{\boldsymbol{\cX}}
\newcommand{\tbcX}{\widetilde{\bcX}}
\newcommand{\tcX}{\widetilde{\cX}}
\newcommand{\cY}{{\mathcal Y}}
\newcommand{\bcY}{\boldsymbol{\cY}}
\newcommand{\cW}{{\mathcal W}}
\newcommand{\bcW}{\boldsymbol{\cW}}
\newcommand{\tcW}{\widetilde{\cW}}

\newcommand{\cc}{\mathrm{c}}
\newcommand{\dd}{\mathrm{d}}
\newcommand{\slu}{{\SC\mathrm{lu}}}

\newcommand{\bbone}{\mathbbm{1}}
\newcommand{\bmid}{\,\big|\,}

\newcommand{\EE}{\operatorname{E}}
\newcommand{\PP}{\operatorname{P}}
\newcommand{\OO}{\operatorname{O}}
\newcommand{\oo}{\operatorname{o}}
\newcommand{\tr}{\operatorname{tr}}
\newcommand{\var}{\operatorname{Var}}
\newcommand{\cov}{\operatorname{Cov}}

\newcommand{\varr}{\varrho}
\newcommand{\vare}{\varepsilon}

\renewcommand{\leq}{\leqslant}
\renewcommand{\geq}{\geqslant}

\newcommand{\stoch}{\stackrel{\PP}{\longrightarrow}}
\newcommand{\distr}{\stackrel{\cL}{\longrightarrow}}
\newcommand{\distre}{\stackrel{\cL}{=}}
\newcommand{\lu}{\stackrel{\slu}{\longrightarrow}}

\newcommand{\nt}{{\lfloor nt\rfloor}}
\newcommand{\nT}{{\lfloor nT\rfloor}}

\newcommand{\proofend}{\hfill\mbox{$\Box$}}

\newcommand{\INARp}{\textup{INAR($p$)}}

\numberwithin{equation}{section}

\newtheorem{Thm}{Theorem}[section]
\newtheorem{Lem}[Thm]{Lemma}

\theorembodyfont{\rm}
\newtheorem{Rem}{Remark}

\begin{document}

\begin{center}

  {\bfseries\Large Asymptotic behavior of critical primitive multi-type 
           branching processes with immigration} \\[5mm]

   {\sc\large M\'arton $\text{Isp\'any}^{*,\diamond}$, \ Gyula $\text{Pap}^{**}$}
\end{center}

\vskip0.2cm

\noindent * University of Debrecen, Faculty of Informatics, Pf.~12, 
     H--4010 Debrecen, Hungary;\\
          ** University of Szeged, Bolyai Institute, Aradi v\'ertan\'uk
                  tere 1, H-6720 Szeged, Hungary

\noindent e--mails: ispany.marton@inf.unideb.hu (M. Isp\'any), 
            papgy@math.u-szeged.hu (G. Pap)

$\diamond$ Corresponding author.

\renewcommand{\thefootnote}{}

\begin{abstract}
Under natural assumptions a Feller type diffusion approximation is derived for
 critical multi-type branching processes with immigration when the offspring
 mean matrix is primitive (in other words, positively regular).
Namely, it is proved that a sequence of appropriately scaled random step
 functions formed from a sequence of critical primitive multi-type branching
 processes with immigration converges weakly towards a squared Bessel process
 supported by a ray determined by the Perron vector of the offspring mean
 matrix.
\end{abstract}

\section{Introduction}
\label{section_intro}

Branching processes have a number of applications in biology, finance,
 economics, queueing theory etc., see e.g.\ Haccou, Jagers and Vatutin
 \cite{HJV}.
Many aspects of applications in epidemiology, genetics and cell kinetics were
 presented at the 2009 Badajoz Workshop on Branching Processes, see
 \cite{Badajoz}.

Let \ $(X_k)_{k \in \ZZ_+}$ \ be a single-type Galton--Watson branching process
 with immigration and with initial value \ $X_0 = 0$.
\ Suppose that it is critial, i.e., the offspring mean equals 1.
Wei and Winnicki \cite{WW} proved a functional limit theorem  
 \ $\cX^{(n)} \distr \cX$ \ as \ $n \to \infty$, \ where
 \ $\cX^{(n)}_t := n^{-1} X_{\nt}$ \ for \ $t \in \RR_+$, \ $n \in \NN$, \ where
 \ $\lfloor x \rfloor$ \ denotes the integer part of \ $x \in \RR$, \ and
 \ $(\cX_t)_{t \in \RR_+}$ \ is a (nonnegative) diffusion process with initial
 value \ $\cX_0 = 0$ \ and with generator
 \begin{equation}\label{Generator}
   Lf(x) = m_\vare f'(x) + \frac{1}{2} V_\xi x f''(x) ,
   \qquad f \in C^\infty_{\cc}(\RR_+) ,
 \end{equation}
 where \ $m_\vare$ \ is the immigration mean, \ $V_\xi$ \ is the offspring
 variance, and \ $C^\infty_{\cc}(\RR_+)$ \ denotes the space of infinitely
 differentiable functions on \ $\RR_+$ \ with compact support.
The process \ $(\cX_t)_{t\in\RR_+}$ \ can also be characterized as the unique
 strong solution of the stochastic differential equation (SDE)
 \[
   \dd \cX_t = m_\vare \, \dd t + \sqrt{ V_\xi \cX_t^+ } \, \dd \cW_t ,
   \qquad t \in \RR_+ ,
 \]
 with initial value \ $\cX_0 = 0$, \ where \ $(\cW_t)_{t\in\RR_+}$ \ is a
 standard Wiener process, and \ $x^+$ \ denotes the positive part of
 \ $x \in \RR$.
\ Note that this so-called square-root process is also known as
 Cox--Ingersoll--Ross model in financial mathematics (see Musiela and Rutkowski
 \cite[p.\ 290]{MR}).
In fact, $(4 V_\xi^{-1} \cX_t)_{t\in\RR_+}$ \ is the square of a
 \ $4 V_\xi^{-1} m_\vare$-dimensional Bessel process started at 0 (see Revuz and
 Yor \cite[XI.1.1]{RY}).

Moreover, for critical Galton--Watson branching processes without immigration,
 Feller \cite{F} proved the following diffusion approximation (see also Ethier
 and Kurtz \cite[Theorem 9.1.3]{EK}).
Consider a sequence of critical Galton--Watson branching processes
 \ $\bigl(X^{(n)}_k\bigr)_{k \in \ZZ_+}$, \ $n \in \NN$, \ without immigration,
 with the same offspring distribution, and with initial value \ $X^{(n)}_0$
 \ independent of the offspring variables such that
 \ $n^{-1} X^{(n)}_0 \distr \mu$ \ as \ $n \to \infty$.
\ Then \ $\cX^{(n)} \distr \cX$ \ as \ $n \to \infty$, \ where
 \ $\cX^{(n)}_t := n^{-1} X^{(n)}_{\nt}$ \ for \ $t \in \RR_+$, \ $n \in \NN$, \ and
 \ $(\cX_t)_{t \in \RR_+}$ \ is a (nonnegative) diffusion process with initial
 distribution \ $\mu$ \ and with generator given by \eqref{Generator} with
 \ $m_\vare = 0$. 

A multi-type branching process \ $(\bX_k)_{k \in \ZZ_+}$ \ is referred to
 respectively as subcritical, critical or supercritical if
 \ $\varr(\bm_\bxi) < 1$, \ $\varr(\bm_\bxi) = 1$ \ or \ $\varr(\bm_\bxi) > 1$,
 \ where \ $\varr(\bm_\bxi)$ \ denotes the spectral radius of the offspring
 mean matrix \ $\bm_\bxi$ \ (see, e.g., Athreya and Ney \cite{AN} or
 Quine \cite{Q}).
Joffe and M\'etivier \cite[Theorem 4.3.1]{JM} studied a sequence
 \ $(\bX^{(n)}_k)_{k \in \ZZ_+}$ \ of critical multi-type branching processes with
 the same offspring distributions but without immigration if the offspring 
 mean matrix is primitive and \ $n^{-1} \bX^{(n)}_0 \distr \bmu$
 \ as \ $n \to \infty$. 
\ They determined the limiting behavior of the martingale part
 \ $(\bcM^{(n)})_{n \in \NN}$ \ given by
 \ $\bcM^{(n)}_t := n^{-1} \sum_{k=1}^{\nt} \bM^{(n)}_k$ \ with
 \ $\bM^{(n)}_k
    := \bX^{(n)}_k - \EE(\bX^{(n)}_k \mid \bX^{(n)}_0, \dots, \bX^{(n)}_{k-1})$
 \ (see \eqref{Conv_M}). Joffe and M\'etivier \cite[Theorem 4.2.2]{JM}
 also studied a sequence \ $(\bX_k^{(n)})_{k \in \ZZ_+}$, \ $n\in\NN$, \ of
 multi-type branching processes without immigration which is nearly critical
 of special type, namely, when the offspring mean matrices \ $\bm_\bxi^{(n)}$, 
 \ $n\in\NN$, \ satisfy \ $\bm_\bxi^{(n)}=\bI_p+n^{-1}\bC+\oo(n^{-1})$ \ as
 \ $n\to\infty$, and they proved that the sequence 
 \ $(n^{-1}\bX_{\nt}^{(n)})_{t \in \RR_+}$ \ converges towards a diffusion process.

The aim of the present paper is to obtain a joint generalization of the above
 mentioned results for critical multi-type branching processes with
 immigration.
We succeeded to determine the asymptotic behavior of a sequence of critical
 multi-type branching processes with immigration and with the same offspring
 and immigration distributions if the offspring mean matrix is primitive and
 \ $n^{-1} \bX^{(n)}_0 \distr \bmu$ \ as \ $n \to \infty$, \ where \ $\bmu$ \ is
 concentrated on the ray \ $\RR_+ \cdot \bu_{\bm_\bxi}$, \ where \ $\bu_{\bm_\bxi}$
 \ is the Perron eigenvector of the offspring mean matrix \ $\bm_\bxi$ \ (see
 Theorem \ref{main}).
It turned out that the limiting diffusion process is always one-dimensional in
 the sense that for all \ $t \in \RR_+$, \ the distribution of \ $\bcX_t$ \ is
 also concentrated on the ray \ $\RR_+ \cdot \bu_{\bm_\bxi}$.
\ In fact, \ $\bcX_t = \cX_t \bu_{\bm_\bxi}$, \ $t \in \RR_+$, \ where
 \ $(\cX_t)_{t \in \RR_+}$ \ is again a squared Bessel process which is
a continuous time and continuous state branching process with immigration.
In the single-type case, Li \cite{Li} proved a result on the convergence of a
sequence of discrete branching processes with immigration to a continuous 
branching process with immigration using appropriate time scaling which is
different from our scaling. Later, Ma \cite{Ma} extended Li's result for 
two-type branching processes. They proved the convergence of the sequence
of infinitesimal generators of single(two)-type branching processes with 
immigration towards the generator of the limiting diffusion process which 
is a well-known technique in case of time-homogeneous Markov processes, 
see, e.g., Ethier and Kurtz \cite{EK}. 
Contrarily, our approach is based on the martingale method.
It is interesting to note that Kesten and Stigum \cite{KS} considered a
 supercritical multi-type branching process without immigration, with a fixed
 initial distribution and with primitive offspring mean matrix, and they
 proved that \ $\varrho(\bm_\xi)^{-n} \bX_n \to \bW$ \ almost surely as
 \ $n \to \infty$, \ where the random vector \ $\bW$ \ is also concentrated on
 the ray \ $\RR_+ \cdot \bu_{\bm_\xi}$
 \ (see also Kurtz, Lyons, Pemantle and Peres \cite{KLPY}). 

\section{Multi-type branching processes with immigration}

Let \ $\ZZ_+$, \ $\NN$, \ $\RR$, \ $\RR_+$  \ and \ $\RR_{++}$ \ denote the set
 of non-negative integers, positive integers, real numbers, non-negative real
 numbers and positive real numbers, respectively.
Every random variable will be defined on a fixed probability space
 \ $(\Omega,\cA,\PP)$.

We will investigate a sequence \ $\bigl(\bX^{(n)}_k\bigr)_{ k \in \ZZ_+ }$,
 \ $n \in \NN$, \ of critical $p$-type branching processes with immigration
 sharing the same offspring and immigration distributions, but having possibly
 different initial distributions.
For each \ $n \in \NN$, \ $k \in \ZZ_+$, \ and \ $i \in \{ 1, \dots, p \}$,
 \ the number of individuals of type \ $i$ \ in the \ $k^\mathrm{th}$
 \ generation of the \ $n^\mathrm{th}$ \ process is denoted by \ $X_{k,i}^{(n)}$.
\ By \ $\xi^{(n)}_{k,j,i,\ell}$ \ we denote the number of type \ $\ell$ \ offspring
 produced by the \ $j^\mathrm{th}$ \ individual who is of type \ $i$
 \ belonging to the \ $(k-1)^\mathrm{th}$ \ generation of the \ $n^\mathrm{th}$
 \ process.
The number of type \ $i$  \ immigrants in the \ $k^\mathrm{th}$ \ generation of
 the \ $n^\mathrm{th}$ \ process will be denoted by \ $\vare^{(n)}_{k,i}$.
\ Consider the random vectors
 \[
   \bX^{(n)}_k := \begin{bmatrix}
                  X_{k,1}^{(n)} \\
                  \vdots \\
                  X_{k,p}^{(n)}
                 \end{bmatrix} , \qquad
   \bxi^{(n)}_{k,j,i} := \begin{bmatrix}
                  \xi^{(n)}_{k,j,i,1} \\
                  \vdots \\
                  \xi^{(n)}_{k,j,i,p}
                 \end{bmatrix} , \qquad
   \bvare^{(n)}_k := \begin{bmatrix}
                \vare^{(n)}_{k,1} \\
                \vdots \\
                \vare^{(n)}_{k,p}
               \end{bmatrix} .
 \]
Then, for \ $n, k \in \NN$, \ we have
 \begin{equation}\label{MBPI(d)}
  \bX^{(n)}_k
  = \sum_{i=1}^p \sum_{j=1}^{X_{k-1,i}^{(n)}} \bxi^{(n)}_{k,j,i} + \bvare^{(n)}_k .
 \end{equation}
Here
 \ $\Bigl\{ \bX^{(n)}_0, \, \bxi^{(n)}_{k,j,i}, \, \bvare^{(n)}_k
            : k, j \in \NN, \, i \in \{ 1, \dots, p \} \Bigr\}$
 \ are supposed to be independent for all \ $n \in \NN$. 
\ Moreover, \ $\big\{ \bxi^{(n)}_{k,j,i} : k, j, n \in \NN \big\}$ \ for each
 \ $i \in \{ 1, \dots, p \}$, \ and \ $\{ \bvare^{(n)}_k : k, n \in \NN \}$
 \ are supposed to consist of identically distributed vectors.

We suppose \ $\EE\bigl(\|\bxi^{(1)}_{1,1,i}\|^2\bigr) < \infty$ \ for all
 \ $i \in \{ 1, \dots, p \}$ \ and
 \ $\EE\bigl(\|\bvare^{(1)}_1\|^2\bigr) < \infty$.
\ Introduce the notations
 \begin{gather*}
  \bm_{\bxi_i} := \EE\bigl(\bxi^{(1)}_{1,1,i}\bigr) \in \RR^p_+ , \qquad
  \bm_{\bxi} := \begin{bmatrix}
                \bm_{\bxi_1} & \cdots & \bm_{\bxi_d}
               \end{bmatrix} \in \RR^{p \times p}_+ , \qquad
  \bm_{\bvare} := \EE\bigl(\bvare^{(1)}_1\bigr) \in \RR^p_+ , \\
  \bV_{\bxi_i} := \var\bigl(\bxi^{(1)}_{1,1,i}\bigr) \in \RR^{p \times p} , \qquad
  \bV_{\bvare} := \var\bigl(\bvare^{(1)}_1\bigr) \in \RR^{p \times p} .
 \end{gather*}
Note that many authors define the offspring mean matrix as \ $\bm^\top_\bxi$. 
\ For \ $k \in \ZZ_+$, \ let
 \ $\cF^{(n)}_k
    := \sigma\bigl( \bX^{(n)}_0, \bX^{(n)}_1 , \dots, \bX^{(n)}_k \bigr)$.
\ By \eqref{MBPI(d)}, 
 \begin{equation}\label{mart}
  \EE\bigl( \bX^{(n)}_k \, \big| \, \cF^{(n)}_{k-1} \bigr) 
  = \sum_{i=1}^p X_{k-1,i}^{(n)} \bm_{\bxi_i} + \bm_{\bvare}
  = \bm_{\bxi} \bX^{(n)}_{k-1} + \bm_{\bvare} .
 \end{equation}
Consequently,
 \begin{equation}\label{recEX}
  \EE\bigl( \bX^{(n)}_k \bigr) 
  = \bm_{\bxi} \EE\bigl( \bX^{(n)}_{k-1} \bigr) + \bm_{\bvare} , \qquad
  k, n \in \NN ,
 \end{equation}
 which implies
 \begin{equation}\label{EXk}
  \EE\bigl( \bX^{(n)}_k \bigr) 
  = \bm_{\bxi}^k \EE\bigl( \bX^{(n)}_0 \bigr)
    + \sum_{j=0}^{k-1} \bm_{\bxi}^j \bm_{\bvare} , \qquad
  k, n \in \NN .
 \end{equation}
Hence, the offspring mean matrix \ $\bm_{\bxi}$ \ plays a crucial role in the
 asymptotic behavior of the sequence \ $\bigl(\bX^{(n)}_k\bigr)_{ k \in \ZZ_+ }$.

In what follows we recall some known facts about primitive nonnegative
 matrices.
A matrix \ $\bA \in \RR^{p \times p}_+$ \ is called primitive if there exists
 \ $m \in \NN$ \ such that \ $\bA^m \in \RR^{p \times p}_{++}$.
\ A matrix \ $\bA \in \RR^{p \times p}_+$ \ is primitive if and only if it is
 irreducible and has only one eigenvalue of maximum modulus; see, e.g., Horn
 and Johnson \cite[Definition 8.5.0, Theorem 8.5.2]{HJ}.
If a matrix \ $\bA \in \RR^{p \times p}_+$ \ is primitive then, by the
 Frobenius--Perron theorem (see, e.g., Horn and Johnson
 \cite[Theorems 8.2.11 and 8.5.1]{HJ}), the following assertions hold:
 \begin{itemize}
  \item
   $\varrho(\bA) \in \RR_{++}$, \ $\varrho(\bA)$ \ is an eigenvalue of \ $\bA$,
    \ the algebraic and geometric multiplicities of \ $\varrho(\bA)$ \ equal 1
    and the absolute values of the other eigenvalues of \ $\bA$ \ are less than
    \ $\varrho(\bA)$.
  \item
   Corresponding to the eigenvalue \ $\varrho(\bA)$ \ there exists a unique
    (right) eigenvector \ $\bu_\bA \in \RR^p_{++}$, \ called Perron vector,
    such that the sum of the coordinates of \ $\bu_\bA$ \ is 1.
  \item
   Further, 
    \[
      \varrho(\bA)^{-n} \bA^n
      \to \bPi_\bA := \bu_\bA \bv_\bA^\top \in \RR^{p \times p}_{++} \qquad
      \text{as \ $n \to \infty$,} 
    \]
    where \ $\bv_\bA \in \RR^p_{++}$ \ is the unique left eigenvector
    corresponding to the eigenvalue \ $\varrho(\bA)$ \ with
    \ $\bu_\bA^\top \bv_\bA = 1$.
  \item
   Moreover, there exist \ $c_\bA , r_\bA \in \RR_{++}$ \ with \ $r_\bA < 1$
    \ such that 
    \begin{equation}\label{rate}
     \| \varrho(\bA)^{-n} \bA^n - \bPi_\bA \| \leq c_\bA r_\bA^n \qquad
     \text{for all \ $n \in \NN$,}
    \end{equation}
    where \ $\|\bBB\|$ \ denotes the operator norm of a matrix
    \ $\bBB \in \RR^{p \times p}$ \ defined by
    \ $\|\bBB\| := \sup_{\|\bx\| = 1} \| \bBB \bx \|$.
 \end{itemize}
A multi-type branching process with immigration will be called primitive if its
 offspring mean matrix \ $\bm_\bxi$ \ is primitive.
Note that many authors call it positively regular.

\section{Convergence results}

A function \ $f : \RR_+ \to \RR^p$ \ is called \emph{c\`adl\`ag} if it is right
 continuous with left limits.
\ Let \ $\DD(\RR_+, \RR^p)$ \ and \ $\CC(\RR_+, \RR^p)$ \ denote the space of
 all $\RR^p$-valued c\`adl\`ag and continuous functions on \ $\RR_+$,
 \ respectively.
Let \ $\cD_\infty(\RR_+, \RR^p)$ \ denote the Borel $\sigma$-algebra in
 \ $\DD(\RR_+, \RR^p)$ \ for the metric defined in Jacod and Shiryaev
 \cite[Chapter VI, (1.26)]{JSh} (with this metric \ $\DD(\RR_+, \RR^p)$ \ is a
 complete and separable metric space).
For $\RR^p$-valued stochastic processes \ $(\bcY_t)_{t \in \RR_+}$ \ and
 \ $(\bcY^{(n)}_t)_{t \in \RR_+}$, \ $n \in \NN$, \ with c\`adl\`ag paths we write
 \ $\bcY^{(n)} \distr \bcY$ \ if the distribution of \ $\bcY^{(n)}$ \ on the
 space \ $(\DD(\RR_+, \RR^p), \cD_\infty(\RR_+, \RR^p))$ \ converges weakly to
 the distribution of \ $\bcY$ \ on the space
 \ $(\DD(\RR_+, \RR^p), \cD_\infty(\RR_+, \RR^p))$ \ as \ $n \to \infty$.

For each \ $n \in \NN$, \ consider the random step processes
 \[
   \bcX^{(n)}_t := n^{-1} \bX^{(n)}_{\nt} ,
   \qquad t \in \RR_+ , \quad n \in \NN .
 \]
For a vector \ $\balpha = (\alpha_i)_{i=1,\dots,p} \in \RR^p_+$, \ we will use
 notation
 \ $\balpha \odot \bV_\bxi
    := \sum_{i=1}^p \alpha_i \bV_{\bxi_i} \in \RR^{p \times p}$,
 \ which is a positive semi-definite matrix, a mixture of the variance
 matrices \ $\bV_{\xi_1}, \ldots, \bV_{\xi_p}$.

\begin{Thm}\label{main}
Let \ $\bigl(\bX^{(n)}_k\bigr)_{ k \in \ZZ_+ }$, \ $n\in\NN$, \ be a sequence of
 critical primitive $p$-type branching processes with immigration sharing the
 same offspring and immigration distributions, but having possibly different
 initial distributions, such that
 \ $n^{-1} \bX^{(n)}_0 \distr \cX_0 \bu_{\bm_\bxi}$, \ where \ $\cX_0$ \ is a
 nonnegative random variable with distribution \ $\mu$.
\ Suppose \ $\EE\bigl(\|\bX^{(n)}_0\|^2\bigr) = \OO(n^2)$,
 \ $\EE\bigl(\|\bxi^{(1)}_{1,1,i}\|^2\bigr) < \infty$ \ for all
 \ $i \in \{ 1, \dots, p \}$ \ and
 \ $\EE\bigl(\|\bvare^{(1)}_1\|^2\bigr) < \infty$.
\ Then
 \begin{gather}\label{Conv_X}
  \bcX^{(n)} \distr \cX \bu_{\bm_\bxi} \qquad \text{as \ $n \to \infty$,}
 \end{gather}
 where \ $(\cX_t)_{t \in \RR_+}$ \ is the unique weak solution (in the sense of
 probability law) of the SDE 
 \begin{equation}\label{SDE_X}
  \dd \cX_t
  = \bv_{\bm_\bxi}^\top \bm_\bvare \, \dd t
    + \sqrt{ \bv_{\bm_\bxi}^\top (\bu_{\bm_\bxi} \odot \bV_\bxi) \bv_{\bm_\bxi}\cX_t^+ }
      \, \dd \cW_t ,
  \qquad t \in \RR_+ ,
 \end{equation}
 with initial distribution \ $\mu$, \ where \ $(\cW_t)_{t \in \RR_+}$ \ is a
 standard Wiener process.
\end{Thm}

\begin{Rem}
We will carry out the proof of Theorem \ref{main} under the assumptions
 \ $\EE\bigl(\|\bxi^{(1)}_{1,1,i}\|^4\bigr) < \infty$ \ for all
 \ $i \in \{ 1, \dots, p \}$ \ and
 \ $\EE\bigl(\|\bvare^{(1)}_1\|^4\bigr) < \infty$.
\ In fact, these higher moment assumptions are needed only for facilitating of
 checking the conditional Lindeberg condition, namely, condition (ii) of
 Theorem \ref{Conv2DiffThm} for proving convergence \eqref{Conv_M} of the
 martingale part. 
One can check the conditional Lindeberg condition under the weaker moment
 assumptions of Theorem \ref{main} by the method of Isp\'any and Pap
 \cite{IP}, see also this method in Barczy et al.~\cite{BarIspPap0}.
If \ $d \geq 2$ \ then it is not clear if one might get rid of the assumption
 \ $\EE(\|\bX^{(n)}_0\|^2) = \OO(n^2)$ \ in Theorem \ref{main}.
\end{Rem}

\begin{Rem}
Under the assumptions of Theorem \ref{main}, by the same method, one can also
 prove \ $\tbcX^{(n)} \distr \tcX \bu_{\bm_\bxi}$ \ as \ $n \to \infty$, \ where
 \ $\tbcX^{(n)}_t
    := n^{-1} \bigl( \bX^{(n)}_{\nt} - \bm_\bxi^{\nt} \bX^{(n)}_0 \bigr)$,
 \ $t \in \RR_+$, \ $n \in \NN$, \  and \ $(\tcX_t)_{t \in \RR_+}$ \ is the unique
 strong solution of the SDE \eqref{SDE_X} with initial value \ $\tcX_0 = 0$.
\end{Rem}

\begin{Rem}\label{Rem_SDE_X}
The SDE \eqref{SDE_X} has a unique strong solution
 \ $(\cX_t^{(x_0)})_{t \in \RR_+}$ \ for all initial values
 \ $\cX_0^{(x_0)} = x_0 \in \RR$.
\ Indeed, since \ $|\sqrt{x} - \sqrt{y}| \leq \sqrt{|x - y|}$ \ for
 \ $x, y \geq 0$, \ the coefficient functions
 \ $\RR \ni x \mapsto \bv_{\bm_\bxi}^\top \bm_\bvare \in \RR_+$ \ and
 \ $\RR \ni x
    \mapsto
    \sqrt{ \bv_{\bm_\bxi}^\top (\bu_\bxi \odot \bV_\bxi) \bv_{\bm_\bxi} x^+ }$
 \ satisfy conditions of part (ii) of Theorem 3.5 in Chapter IX in Revuz and
 Yor \cite{RY} or the conditions of Proposition 5.2.13 in Karatzas and
 Shreve \cite{KarShr}.
Further, by the comparison theorem
 (see, e.g., Revuz and Yor \cite[Theorem 3.7, Chapter IX]{RY}), if the
 initial value \ $\cX_0^{(x_0)} = x_0$ \ is nonnegative, then \ $\cX_t^{(x)}$ \ is
 nonnegative for all \ $t \in \RR_+$ \ with probability one.
Hence \ $\cX_t^+$ \ may be replaced by \ $\cX_t$ \ under the square root in
 \eqref{SDE_X}.
\end{Rem}

\noindent
\textit{Proof of Theorem \ref{main}.} \
In order to prove \eqref{Conv_X}, for each \ $n \in \NN$, \ introduce the
 sequence
 \begin{equation}\label{Mk}
  \bM^{(n)}_k := \bX^{(n)}_k - \EE\bigl(\bX^{(n)}_k \, \big| \, \cF^{(n)}_{k-1}\bigr)
             = \bX^{(n)}_k - \bm_{\bxi} \bX^{(n)}_{k-1} - \bm_{\bvare} ,
  \qquad k \in \NN ,
 \end{equation}
 which is a sequence of martingale differences with respect to the filtration
 \ $\bigl(\cF^{(n)}_k\bigr)_{k \in \ZZ_+}$.
\ Consider the random step processes
 \[
   \bcM^{(n)}_t := n^{-1} \bigg( \bX^{(n)}_0 + \sum_{k=1}^{\nt} \bM^{(n)}_k \bigg) ,
   \qquad t \in \RR_+ , \quad n \in \NN .
 \]
First we will verify convergence
 \begin{gather}\label{Conv_M}
  \bcM^{(n)} \distr \bcM \qquad \text{as \ $n \to \infty$,}
 \end{gather}
 where \ $(\bcM_t)_{t \in \RR_+}$ \ is the unique weak solution of the SDE
 \begin{equation}\label{SDE_M}
  \dd \bcM_t
  = \sqrt{ (\bPi_{\bm_\bxi} (\bcM_t + t \bm_\bvare))^+ \odot \bV_\bxi } \,
    \dd \bcW_t ,
  \qquad t \in \RR_+ ,
 \end{equation}
 with initial distribution \ $\bmu :\distre \cX_0 \bu_{\bm_\bxi}$, \ where
 \ $(\bcW_t)_{t\in\RR_+}$ \ is a standard $p$-dimensional Wiener process, 
 \ $\bx^+$ \ denotes the positive part of \ $\bx \in \RR^p$, \ and for a
 positive semi-definite matrix \ $\bA \in \RR^{p \times p}$, \ $\sqrt{\bA}$
 \ denotes its unique symmetric positive semi-definite square root.

From \eqref{Mk} we obtain the recursion
 \begin{equation}\label{regr}
  \bX^{(n)}_k =  \bm_{\bxi} \bX^{(n)}_{k-1} + \bM^{(n)}_k + \bm_{\bvare} ,
  \qquad k \in \NN ,
 \end{equation} 
 implying
 \begin{equation}\label{X}
  \bX^{(n)}_k = \bm_{\bxi}^k \bX^{(n)}_0
              + \sum_{j=1}^k \bm_{\bxi}^{k-j} (\bM^{(n)}_j + \bm_{\bvare}) ,
  \qquad k \in \NN .
 \end{equation} 
Applying a version of the continuous mapping theorem (see Appendix) together
 with \eqref{Conv_M} and \eqref{X}, in Section \ref{Proof} we show that
 \begin{gather}\label{Conv_bX}
  \bcX^{(n)} \distr \bcX \qquad \text{as \ $n \to \infty$,}
 \end{gather}
 where \ $\bcX_t := \bPi_{\bm_\bxi} (\bcM_t + t \bm_\bvare)$, \ $t \in \RR_+$.
\ Using \ $\bPi_{\bm_\bxi} = \bu_{\bm_\bxi} \bv_{\bm_\bxi}^\top$ \ and
 \ $\bv_{\bm_\bxi}^\top \bu_{\bm_\bxi} = 1$ \ we get that the process
 \ $\cY_t := \bv_{\bm_\bxi}^\top \bcX_t$, \ $t \in \RR_+$, \ satisfies
 \ $\cY_t = \bv_{\bm_\bxi}^\top \bPi_{\bm_\bxi} (\bcM_t + t \bm_\bvare)
    = \bv_{\bm_\bxi}^\top (\bcM_t + t \bm_\bvare)$, \ $t \in \RR_+$,
 \ hence \ $\bcX_t = \cY_t \bu_{\bm_\bxi}$.
\ By It\^o's formula we obtain that \ $(\cY_t)_{t \in \RR_+}$ \ satisfies the SDE
 \eqref{SDE_X} (see the analysis of the process \ $(\cP_t^{(\by_0)})_{t \in \RR_+}$
 \ in the first equation of \eqref{SDE_P_Q} and in equation \eqref{CIR3}) such
 that
 \ $\cY_0 = \bv_{\bm_\bxi}^\top \bcX_0 = \bv_{\bm_\bxi}^\top (\cX_0 \bu_{\bm_\bxi})
          = \cX_0$,
 \ thus we conclude the statement of Theorem \ref{main}. 
\proofend

\begin{Rem}
By It\^o's formula, the limit process \ $(\bcX_t)_{t \in \RR_+}$ \ in
 \eqref{Conv_X} can also be characterized as a weak solution of the SDE 
 \begin{equation}\label{SDE_bX}
  \dd \bcX_t
  = \bPi_{\bm_\bxi} \bm_\bvare \, \dd t
    + \bPi_{\bm_\bxi} \sqrt{ \bcX_t^+ \odot \bV_\bxi } \,
      \dd \bcW_t ,
  \qquad t \in \RR_+ ,
 \end{equation}
 with initial distribution
 \ $\bPi_{\bm_\bxi} \bcM_0 = \bPi_{\bm_\bxi} (\cX_0 \bu_{\bm_\bxi})
    = \cX_0 \bu_{\bm_\bxi}$,
 \ since 
 \ $\bPi_{\bm_\bxi} \bu_{\bm_\bxi} = \bu_{\bm_\bxi} \bv_{\bm_\bxi}^\top \bu_{\bm_\bxi}
    = \bu_{\bm_\bxi}$.
\end{Rem}

\begin{Rem}
The generator of \ $(\bcM_t)_{t \in \RR_+}$ \ is given by
 \begin{align*}
  L_t f(\bx)
  &= \frac{1}{2}
     \big\langle [(\bPi_{\bm_\bxi} (\bx + t \bm_\bvare)) \odot \bV_\bxi] \nabla ,
                   \nabla \big\rangle
     f(\bx) \\
  &= \frac{1}{2} (\bx + t \bm_\bvare)^\top \bPi_{\bm_\bxi}^\top
     \sum_{i=1}^p \sum_{j=1}^p
      \bV_{\xi,i,j} \, \partial_i \partial_j f(\bx) ,
  \qquad t \in \RR_+ , \quad f \in C^\infty_{\cc}(\RR^p) ,
 \end{align*}
 where
 \ $\bV_{\xi,i,j} := (\cov(\xi_{1,1,\ell,i}, \xi_{1,1,\ell,j}))_{\ell=1,\ldots,d}
                \in \RR^p_+$. 
\ (Joffe and M\'etivier \cite[Theorem 4.3.1]{JM} also obtained this generator
 with \ $\bm_\bvare = \bzero$ \ deriving \eqref{Conv_M} for processes without
 immigration.)
\end{Rem}

\section{Proof of \ $\bcM^{(n)} \distr \bcM$ \ and
            \ $\bcX^{(n)} \distr \bcX$ \ as \ $n \to \infty$}
\label{Proof}

First we prove \ $\bcM^{(n)} \distr \bcM$ \ applying Theorem \ref{Conv2DiffThm}
 for \ $\bcU=\bcM$, \ $\bU^{(n)}_0 = n^{-1} \bX^{(n)}_0$ \ and
 \ $\bU^{(n)}_k = n^{-1} \bM^{(n)}_k$ \ for \ $n, k \in \NN$, \ and with
 coefficient function \ $\gamma : \RR_+ \times \RR^p \to \RR^{p \times p}$ \ of
 the SDE \eqref{SDE_M} given by
 \ $\gamma(t, \bx)
    = \sqrt{ (\bPi_{\bm_\bxi} (\bx + t \bm_\bvare))^+ \odot \bV_\bxi}$.
\ 
The aim of the following discussion is to show that the SDE \eqref{SDE_M} has
 a unique strong solution \ $\bigl(\bcM_t^{(\by_0)}\bigr)_{t \in \RR_+}$ \ with
 initial value \ $\bcM_0^{(\by_0)} = \by_0$ \ for all \ $\by_0 \in \RR^p$. 
\ First suppose that the SDE \eqref{SDE_M}, which can also be written in the
 form
 \[
   \dd \bcM_t
  = \sqrt{ (\bv_{\bm_\bxi}^\top (\bcM_t + t \bm_\bvare))^+
           (\bu_{\bm_\bxi} \odot \bV_\bxi) } \,
    \dd \bcW_t ,
 \]
 has a strong solution \ $\bigl(\bcM_t^{(\by_0)}\bigr)_{t \in \RR_+}$ \ with
 \ $\bcM_0^{(\by_0)} = \by_0$. 
\ Then, by It\^o's formula, the process
 \ $\bigl(\cP_t^{(\by_0)}, \, \bcQ_t^{(\by_0)}\bigr)_{t \in \RR_+}$, \ defined by
 \[
   \cP_t^{(\by_0)} := \bv_{\bm_\bxi}^\top (\bcM_t^{(\by_0)} + t \bm_\bvare) , \qquad
   \bcQ_t^{(\by_0)} := \bcM_t^{(\by_0)} - \cP_t^{(\by_0)} \bu_{\bm_\bxi}
 \]
 is a strong solution of the SDE
 \begin{equation}\label{SDE_P_Q}
  \begin{cases}
   \dd \cP_t
   = \bv_{\bm_{\bxi}}^\top \bm_{\bvare} \, \dd t
     + \sqrt{\cP_t^+} \, \bv_{\bm_{\bxi}}^\top
       \sqrt{\bu_{\bm_{\bxi}} \odot \bV_\bxi} 
       \, \dd \bcW_t, \\[2mm]
   \dd \bcQ_t
   = - \bPi_{\bm_{\bxi}} \bm_{\bvare} \, \dd t
     + \sqrt{\cP_t^+} \, \bigl( \bI_p - \bPi_{\bm_{\bxi}} \bigr)
       \sqrt{\bu_{\bm_{\bxi}} \odot \bV_\bxi} 
       \, \dd \bcW_t
  \end{cases}
 \end{equation} 
 with initial value
 \ $\bigl(\cP_0^{(\by_0)}, \, \bcQ_0^{(\by_0)}\bigr)
    =\big(\bv_{\bm_\bxi}^\top \by_0, \, (\bI_d - \bPi_{\bm_\bxi})\by_0\bigr)$,
 \ where \ $\bI_p$ \ denotes the $p$-dimensional unit matrix.
The SDE \eqref{SDE_P_Q} has a unique strong solution
 \ $\bigl(\cP_t^{(p_0)}, \, \bcQ_t^{(\bq_0)}\bigr)_{t \in \RR_+}$, \ with an
 arbitrary initial value
 \ $\big(\cP_0^{(p_0)}, \, \bcQ_0^{(\bq_0)}\big)
    = (p_0, \bq_0) \in \RR_+ \times \RR^p$, \ since
 the first equation of \eqref{SDE_P_Q} can be written in the form 
 \begin{equation}\label{CIR3}
  \dd \cP_t = b \, \dd t + \sqrt{\cP_t^+} \, \dd \tcW_t
 \end{equation}
 with \ $b := \bv_{\bm_\bxi}^\top \bm_\bvare \in \RR_+$ \ and
 \[
   \tcW_t
   := \bv_{\bm_{\bxi}}^\top \sqrt{\bu_{m_{\xi}} \odot \bV_\bxi} \, \bcW_t
   = \sqrt{\bv_{\bm_\bxi}^\top (\bu_{\bm_\bxi} \odot \bV_\bxi) \bv_{\bm_\bxi}} \,
     \cW_t ,
 \]
 where \ $(\cW_t)_{t \in \RR_+}$ \ is a standard one-dimensional Wiener process.
(Equation \eqref{CIR3} can be discussed as equation \eqref{SDE_X} in Remark
 \eqref{Rem_SDE_X}.)
If \ $(\cP_t^{(\by_0)}, \, \bcQ_t^{(\by_0)})_{t \in \RR_+}$ \ is the unique strong
 solution of the SDE \eqref{SDE_P_Q} with the initial value
 \ $\bigl(\cP_0^{(\by_0)}, \, \bcQ_0^{(\by_0)}\bigr)
    = \bigl(\bv_{\bm_\bxi}^\top \by_0, \, (\bI_p - \bPi_{\bm_\bxi})\by_0\bigr)$,
 \ then, again by It\^o's formula,
 \[
   \bcM_t^{(\by_0)} := \cP_t^{(\by_0)} \, \bu_{\bm_{\bxi}} + \bcQ_t^{(\by_0)} , \qquad
   t \in \RR_+ , 
 \] 
 is a strong solution of \eqref{SDE_M} with \ $\bcM_0^{(\by_0)} = \by_0$.
\ Consequently, \eqref{SDE_M} admits a unique strong solution
 \ $\bigl(\bcM_t^{(\by_0)}\bigr)_{t \in \RR_+}$ \ with \ $\bcM_0^{(\by_0)} = \by_0$
 \ for all \ $\by_0 \in \RR^p$.

Now we show that conditions (i) and (ii) of Theorem \ref{Conv2DiffThm} hold.
We have to check that, for each \ $T > 0$,
 \begin{align} \label{Cond1}
  &\sup_{t\in[0,T]}
    \bigg\| \frac{1}{n^2} \sum_{k=1}^{\nt}
            \EE\bigl[ \bM^{(n)}_k (\bM^{(n)}_k)^\top \,
                      \big| \, \cF^{(n)}_{k-1} \bigr]
            - \int_0^t (\bcR^{(n)}_s)_+\, \dd s \, \odot \bV_\bxi \bigg\|
   \stoch 0,\\
  &\frac{1}{n^2}
   \sum_{k=1}^{\nT}
    \EE\bigl( \|\bM^{(n)}_k\|^2 \bbone_{\{\|\bM^{(n)}_k\| > n \theta\}} \,
              \big| \, \cF^{(n)}_{k-1} \bigr)
   \stoch 0   \qquad\text{for all \ $\theta>0$} \label{Cond2}
 \end{align}
 as \ $n \to \infty$, \ where the process \ $(\bcR^{(n)}_t)_{t \in \RR_+}$ \ is
 defined by
 \begin{equation}\label{Rnt}
  \bcR^{(n)}_t := \bPi_{\bm_\bxi} \bigl(\bcM^{(n)}_t + t \bm_\bvare\bigr) , \qquad
  t \in \RR_+ , \quad n \in \NN .
 \end{equation}
By \eqref{Mk},
 \begin{align*}
  \bcR^{(n)}_t
  & = \bPi_{m_\xi}
      \bigg( n^{-1}
             \bigg( \bX^{(n)}_0
                    + \sum_{k=1}^{\nt}
                       ( \bX^{(n)}_k - \bm_\bxi \bX^{(n)}_{k-1} - \bm_\bvare )
             \bigg)
             + t \bm_\bvare \bigg) \\
  & = n^{-1} \bPi_{\bm_\bxi} \bX^{(n)}_{\nt}
      + n^{-1} ( nt - \nt ) \bPi_{\bm_\bxi} \bm_\bvare ,
 \end{align*}
 where we used that 
 \[
   \bPi_{\bm_\bxi} \bm_\bxi
   = \left( \lim\limits_{n \to \infty} \bm_\bxi^n \right) \bm_\bxi
   = \lim\limits_{n \to \infty} \bm_\bxi^{n+1} = \bPi_{\bm_\bxi}
 \]
 implies \ $\bPi_{\bm_\bxi} (\bI_p - \bm_\bxi) = \bzero$.
\ Thus \ $(\bcR^{(n)}_t)_+ = \bcR^{(n)}_t$, \ and 
 \begin{align*}
  \int_0^t (\bcR^{(n)}_s)_+ \,\dd s
  = \frac{1}{n^2} \sum_{\ell=0}^{\nt-1} \bPi_{\bm_\bxi} \bX^{(n)}_\ell
    + \frac{nt-\nt}{n^2} \bPi_{\bm_\bxi} \bX^{(n)}_{\nt} 
    + \frac{\nt+(nt-\nt)^2}{2n^2} \bPi_{\bm_\bxi} \bm_\bvare .
 \end{align*}
Using \eqref{Mcond}, we obtain
 \[
   \frac{1}{n^2} \sum_{k=1}^{\nt}
   \EE\bigl[ \bM^{(n)}_k (\bM^{(n)}_k)^\top \, \big| \, \cF^{(n)}_{k-1} \bigr]
   = \frac{\nt}{n^2} \bV_\bvare
     + \frac{1}{n^2} \sum_{k=1}^{\nt} \bX_{k-1}^{(n)} \odot \bV_\bxi .
 \]
Hence, in order to show \eqref{Cond1}, it suffices to prove
 \begin{equation}\label{Cond11}
  \sup_{t \in [0,T]}
   \frac{1}{n^2}
   \sum_{k=0}^{\nt-1}
    \| (\bI_p - \bPi_{\bm_{\bxi}}) \bX^{(n)}_k \|
  \stoch 0,\qquad
  \sup_{t \in [0,T]}
   \frac{1}{n^2}
   \|\bX^{(n)}_{\nt}\|
  \stoch 0
 \end{equation}
 as \ $n \to \infty$.
\ Using \eqref{X} and \ $\bPi_{\bm_\xi} \bm_\bxi = \bPi_{\bm_\bxi}$, \ we obtain
 \[
   (\bI_d - \bPi_{\bm_\bxi}) \bX^{(n)}_k
   = \left( \bm_\bxi^k - \bPi_{\bm_\bxi} \right) \bX^{(n)}_0
     + \sum_{j=1}^k
        \left( \bm_\bxi^{k-j} - \bPi_{\bm_\bxi} \right)
        (\bM^{n)}_j + \bm_\bvare) .
 \]
Hence by \eqref{rate},
 \begin{align*} 
  \sum_{k=0}^{\nt-1}
   \| (\bI_d - \bPi_{\bm_{\bxi}}) \bX^{(n)}_k \|
   & \leq c_{\bm_{\bxi}}
          \sum_{k=0}^{\nt-1}
           r_{\bm_\bxi}^k
           \|\bX^{(n)}_0\|
          + c_{\bm_{\bxi}}
            \sum_{k=1}^{\nt-1}
             \sum_{j=1}^k
              r_{\bm_\bxi}^{k-j}
              \|\bM^{(n)}_j + \bm_\bvare\| \\
   & \leq \frac{c_{\bm_\bxi}}{1 - r_{\bm_\bxi}}
          \bigg( \|\bX^{(n)}_0\| + \nt \cdot \|\bm_\bvare\|
                 + \sum_{j=1}^{\nt-1} \|\bM^{(n)}_j\| \bigg) .
 \end{align*}
Moreover, by \eqref{X} and \eqref{C},
 \begin{align*}
  \|\bX_{\nt}\|
  & \leq \| \bm_\bxi^{\nt} \| \cdot \| \bX^{(n)}_0 \|
         + \sum_{j=1}^{\nt}
            \| \bm_\bxi^{\nt-j} \| \cdot \|\bM^{(n)}_j + \bm_\bvare\| \\
  & \leq C_{\bm_{\bxi}}
         \bigg( \| \bX^{(n)}_0 \| + \nt \cdot \|\bm_\bvare\|
                + \sum_{j=1}^{\nt} \|\bM^{(n)}_j\| \bigg) ,
 \end{align*}
 where \ $C_{\bm_{\bxi}}$ \ is defined by \eqref{C}.
Consequently, in order to prove \eqref{Cond11}, it suffices to show
 \[
  \frac{1}{n^2} \sum_{j=1}^{\nT} \|\bM^{(n)}_j\| \stoch 0, \qquad
  \frac{1}{n^2} \|\bX^{(n)}_0\| \stoch 0 \qquad
  \text{as \ $n \to \infty$.}
 \]
In fact, assumption \ $n^{-1} \bX^{(n)}_0 \distr \bmu$ \ implies the second
 convergence, while Lemma \ref{EEX} yields
 \ $n^{-2} \sum_{j=1}^{\nT} \EE(\|\bM^{(n)}_j\|) \to 0$, \ thus we obtain
 \eqref{Cond1}.

Next we check condition \eqref{Cond2}.
We have
 \[
   \EE\bigl( \|\bM^{(n)}_k\|^2 \bbone_{\{\|\bM^{(n)}_k\| > n \theta\}} \,
              \big| \, \cF^{(n)}_{k-1} \bigr)
   \leq n^{-2} \theta^{-2}
        \EE\bigl( \|\bM^{(n)}_k\|^4 \, \big| \, \cF^{(n)}_{k-1} \bigr) .
 \]  
Moreover, \ $n^{-4} \sum_{k=1}^{\nT} \EE\bigl(\|\bM^{(n)}_k\|^4\bigr) \to 0$ \ as
 \ $n \to \infty$, \ since
 \ $\EE\bigl(\|\bM^{(n)}_k\|^4\bigr) = \OO\bigl((k+n)^2\bigr)$ \ by Lemma
 \ref{EEX}.
Hence we obtain \eqref{Cond2}.

Now we turn to prove \eqref{Conv_bX} applying Lemma \ref{Conv2Funct}.
By \eqref{X}, \ $\bcX^{(n)} = \Psi_n(\bcM^{(n)})$, \ where the mapping
 \ $\Psi_n : \DD(\RR_+, \RR^p) \to \DD(\RR_+, \RR^p)$ \ is given by
 \[
   \Psi_n(f)(t)
   := \bm_\bxi^\nt f(0)
      + \sum_{j=1}^\nt
         \bm_\bxi^{\nt - j}
          \left( f\left(\frac{j}{n}\right) - f\left(\frac{j-1}{n}\right)
                 + n^{-1} \bm_\bvare \right)
 \]
 for \ $f \in \DD(\RR_+, \RR^p)$, \ $t \in \RR_+$, \ $n \in \NN$.
\ Further, \ $\bcX = \Psi(\bcM)$, \ where the mapping
 \ $\Psi : \DD(\RR_+, \RR^p) \to \DD(\RR_+, \RR^p)$ \ is given by
 \[
   \Psi(f)(t) := \bPi_{\bm_\bxi} (f(t) + t \bm_\bvare) , \qquad
   f \in \DD(\RR_+, \RR^p), \qquad t \in \RR_+ .
 \]
Measurability of the mappings \ $\Psi_n$, \ $n \in \NN$, \ and \ $\Psi$ \ can
 be checked as in Barczy et al. \cite{BarIspPap0}.

The aim of the following discussion is to show that the set
 \ $C := \{ f \in \CC(\RR_+, \RR^p) : \bPi_{\bm_\bxi} f(0) = f(0) \}$ \ satisfies
 \ $C \in \cD_\infty(\RR_+, \RR^p)$, \ $C \subset C_{\Psi, \, (\Psi_n)_{n \in \NN}}$
 \ and \ $\PP(\bcM \in C) = 1$.
\ Note that \ $f \in C$ \ implies \ $f(0) \in \RR \cdot \bu_{\bm_\bxi}$.

First note that
 \ $C = \CC(\RR_+, \RR^p)
        \cap \pi_0^{-1}\bigl((\bI_p - \bPi_{\bm_\bxi})^{-1}(\{ \bzero \})\bigr)$,
 \ where \ $\pi_0 : \DD(\RR_+, \RR^p) \to \RR^p$ \ denotes the projection
 defined by \ $\pi_0(f) := f(0)$ \ for \ $f \in \DD(\RR_+, \RR^p)$.
\ Using that \ $\CC(\RR_+, \RR^p) \in \cD_\infty$ \ (see, e.g., Ethier and Kurtz
 \cite[Problem 3.11.25]{EK}), the mapping
 \ $\RR^p \ni \bx \mapsto (\bI_p - \bPi_{\bm_\bxi}) \bx \in \RR^p$ \ is
 measurable and that \ $\pi_0$ \ is measurable (see, e.g., Ethier and Kurtz
 \cite[Proposition 3.7.1]{EK}), we obtain \ $C \in \cD_\infty(\RR_+, \RR^p)$.
 
Fix a function \ $f \in C$ \ and a sequence \ $(f_n)_{n\in\NN}$ \ in
 \ $\DD(\RR^p)$ \ with \ $f_n \lu f$.
\ By the definition of \ $\Psi$, \ we have \ $\Psi(f) \in \CC(\RR^p)$.
\ Further, we can write
 \begin{align*}
  \Psi_n(f_n)(t)
  &= \bPi_{\bm_\bxi}
     \left( f_n\left(\frac{\nt}{n}\right)
            + \frac{\nt}{n} \bm_\bvare \right)
     + \bigl( \bm_\bxi^{\nt} - \bPi_{\bm_\bxi} \bigr) f(0) \\
  &\quad
     + \sum_{j=1}^{\nt}
        \bigl( \bm_\bxi^{\nt-j} - \bPi_{\bm_\bxi} \bigr)
        \left( f_n\left(\frac{j}{n}\right)
               - f\left(\frac{j-1}{n}\right)
               + \frac{1}{n} \bm_\bvare \right) ,
 \end{align*}
 hence we have
 \begin{align*}
  \|\Psi_n(f_n)(t)-\Psi(f)(t)\|
   &\leq \| \bPi_{\bm_\bxi} \|
         \left( \left\| f_n\left(\frac{\nt}{n}\right)
                        - f(t) \right\|
                + \frac{1}{n} \|\bm_\bvare\| \right)
         + \left\| \left( \bm_\bxi^{\nt} - \bPi_{\bm_\bxi} \right)
                   f_n(0) \right\| \\
   &\quad
         + \sum_{j=1}^{\nt}
            \bigl\| \bm_\bxi^{\nt-j} - \bPi_{\bm_\bxi} \bigr\|
            \left( \left\| f_n\left(\frac{j}{n}\right)
                           - f_n\left(\frac{j-1}{n}\right) \right\|
                   + \frac{1}{n} \|\bm_\bvare\| \right) .
 \end{align*}
For all \ $T > 0$ \ and \ $t \in [0, T]$,
 \begin{align*}
  \left\| f_n\left(\frac{\nt}{n}\right) - f(t) \right \|
  &\leq \left\| f_n\left(\frac{\nt}{n}\right)
                - f\left(\frac{\nt}{n}\right) \right\|
        + \left\| f\left(\frac{\nt}{n}\right) - f(t) \right\| \\
  &\leq \omega_T(f, n^{-1}) + \sup_{t \in [0,T]} \|f_n(t) - f(t)\| ,
 \end{align*}
 where \ $\omega_T(f, \cdot)$ \ is the modulus of continuity of \ $f$ \ on
 \ $[0, T]$, \ and we have \ $\omega_T(f, n^{-1}) \to 0$ \ since \ $f$ \ is
 continuous (see, e.g., Jacod and Shiryaev \cite[VI.1.6]{JSh}).
In a similar way,
 \[
   \left\| f_n\left(\frac{j}{n}\right)
           - f_n\left(\frac{j-1}{n}\right) \right\|
   \leq \omega_T(f, n^{-1}) + 2 \sup_{t \in [0,T]} \|f_n(t) - f(t)\| .
 \]
By \eqref{rate},
 \[
   \sum_{j=1}^{\nt} \bigl\| \bm_\bxi^{\nt-j} - \bPi_{\bm_\bxi} \bigr\|
   \leq \sum_{j=1}^{\nT} c_{\bm_\bxi} r_{\bm_\bxi}^{\nt-j}
   \leq \frac{c_{\bm_\bxi}}{1 - r_{\bm_\bxi}} .
 \]
Further,
 \begin{align*}
  \Big\| \big( \bm_\bxi^{\nt} - \bPi_{\bm_\bxi} \big)
           f_n(0) \Big\| 
  & \leq \Big\| \big( \bm_\bxi^{\nt} - \bPi_{\bm_\bxi} \big)
                \big( f_n(0) - f(0) \big) \Big\|
         + \Big\| \big( \bm_\bxi^{\nt} - \bPi_{\bm_\bxi} \big) f(0) \Big\| \\
  & \leq c_{\bm_\bxi} \sup_{t \in [0,T]} \|f_n(t) - f(t)\| ,
 \end{align*}
 since \ $\big( \bm_\bxi^{\nt} - \bPi_{\bm_\bxi} \big) f(0) = \bzero$ \ for all
 \ $t \in \RR_+$.
\ Indeed,
 \ $\bm_\bxi \bPi_{\bm_\bxi} = \bm_\bxi \lim_{n \to \infty} \bm_\bxi^n
    = \lim_{n \to \infty} \bm_\bxi^{n+1} = \bPi_{\bm_\bxi}$ \ and
 \ $f(0) = \bPi_{\bm_\bxi} f(0)$ \ imply
 \ $\bm_\bxi^{\nt} f(0) = \bm_\bxi^{\nt} \bPi_{\bm_\bxi} f(0) = \bPi_{\bm_\bxi} f(0)$.
\ Thus we conclude \ $C \subset C_{\Psi, \, (\Psi_n)_{n \in \NN}}$.

By the definition of a weak solution (see, e.g., Jacod and Shiryaev
 \cite[Definition 2.24, Chapter III]{JSh}), \ $\bcM$ \ has almost sure
 continuous sample paths, so we have \ $\PP(\bcM \in C) = 1$.
\ Consequently, by Lemma \ref{Conv2Funct}, we obtain
 \ $\bcX^{(n)} = \Psi_n(\bcM^{(n)}) \distr \Psi(\bcM) \distre \bcX$ \ as
 \ $n \to \infty$.
\proofend

\section*{Acknowledgements}

The authors have been supported by the Hungarian
 Chinese Intergovernmental S \& T Cooperation Programme for 2011-2013
 under Grant No.\ 10-1-2011-0079.
G. Pap has been partially supported by the Hungarian Scientific Research Fund
 under Grant No.\ OTKA T-079128.
M. Isp\'any has been partially supported by
 T\'AMOP 4.2.1./B-09/1/KONV-2010-0007/IK/IT project, which is implemented
 through the New Hungary Development Plan co-financed by the European Social
 Fund and the European Regional Development Fund.

\section*{Appendix}
\label{app}
\setcounter{section}{1}
\setcounter{equation}{0}
\renewcommand{\thesection}{\Alph{section}}

In the proof of Theorem \ref{main} we will use some facts about the first and
 second order moments of the sequences \ $(\bX^{(n)}_k)_{k \in \ZZ_+}$ \ and
 \ $(\bM^{(n)}_k)_{k \in \NN}$.

\begin{Lem}\label{Moments}
Under the assumptions of Theorem \ref{main} we have for all \ $k, n \in \NN$
 \begin{gather}
  \EE(\bX^{(n)}_k)
  = \bm_{\bxi}^k \EE(\bX^{(n)}_0)
    + \sum_{j=0}^{k-1} \bm_{\bxi}^j \bm_{\bvare}, \label{mean}\\
  \begin{split}
   \var(\bX^{(n)}_k)
   &= \sum_{j=0}^{k-1}
       \bm_{\bxi}^j
       \left[ \bV_{\bvare}
              + (\bm_{\bxi}^{k-j-1} \EE(\bX^{(n)}_0)) \odot \bV_\xi \right]
       (\bm_{\bxi}^\top)^j \\
   & \quad
      + \bm_{\bxi}^k (\var(\bX^{(n)}_0)) (\bm_{\bxi}^\top)^k
      + \sum_{j=0}^{k-2}
         \bm_{\bxi}^j
         \sum_{\ell=0}^{k-j-2}
          \left[ (\bm_{\bxi}^\ell \bm_{\bvare}) \odot \bV_\xi \right]
          (\bm_{\bxi}^\top)^j .
  \end{split} \label{var}     
 \end{gather}
Moreover,
 \begin{align} 
  \EE \bigl( \bM^{(n)}_k \, \big| \, \cF^{(n)}_{k-1} \bigr)
  & = \bzero \qquad \text{for \ $k, n \in \NN$,}
   \label{m} \\ 
  \EE \bigl[ \bM^{(n)}_k (\bM^{(n)}_\ell)^\top
             \, \big| \, \cF^{(n)}_{\max\{k,\ell\}-1} \bigr]
  & = \begin{cases}
       \bV_{\bvare} + \bX_{k-1}^{(n)} \odot \bV_{\bxi}
        & \text{if \ $k = \ell$,} \\
       \bzero & \text{if \ $k \ne \ell$.}
      \end{cases} \label{Mcond}
 \end{align}
Further,
 \begin{align}
  \EE(\bM^{(n)}_k) & = \bzero \qquad \text{for \ $k \in \NN$,} \label{M} \\
  \EE\bigl[ \bM^{(n)}_k (\bM^{(n)}_\ell)^\top \bigr]
   & = \begin{cases}
        \bV_{\bvare} + \EE(\bX_{k-1}^{(n)}) \odot \bV_{\bxi}
         & \text{if \ $k = \ell$,} \\
        \bzero & \text{if \ $k \ne \ell$.}
       \end{cases} \label{Cov}
 \end{align}
\end{Lem}

\noindent
\textit{Proof.}
We have already proved \eqref{mean}, see \eqref{EXk}.
The equality
 \ $\bM^{(n)}_k
    = \bX^{(n)}_k - \EE\bigl(\bX^{(n)}_k \, \big| \, \cF^{(n)}_{k-1}\bigr)$
 \ clearly implies \eqref{m} and \eqref{M}.
By \eqref{MBPI(d)} and \eqref{Mk},
 \begin{equation}\label{Mdeco}
   \bM^{(n)}_k
   = \bX^{(n)}_k - \sum_{i=1}^p X_{k-1,i}^{(n)} \, \EE(\bxi^{(1)}_{1,1,i}) - \bm_{\bvare}
   = ( \bvare_k - \EE(\bvare_k) )
     + \sum_{i=1}^p \sum_{j=1}^{X_{k-1,i}^{(n)}}
        ( \bxi^{(n)}_{k,j,i} - \EE(\bxi^{(n)}_{k,j,i}) ) .
 \end{equation}
For each \ $k, n \in \NN$, \ the random vectors
 \ $\big\{\bxi^{(n)}_{k,j,i} - \EE(\bxi^{(n)}_{k,j,i}) , \,
          \bvare^{(n)}_k - \EE(\bvare^{(n)}_k)
          : j \in \NN , \, i \in \{1, \dots, p\} \big\}$
 \ are independent of each others, independent of \ $\cF^{(n)}_{k-1}$, \ and have
 zero mean, thus in case \ $k = \ell$ \ we conclude \eqref{Mcond} and hence
 \eqref{Cov}.
If \ $k < \ell$ \ then
 \ $\EE\bigl[ \bM^{(n)}_k (\bM^{(n)}_\ell)^\top \, \big| \, \cF^{(n)}_{\ell-1} \bigr]
    = \bM^{(n)}_k
      \EE\bigl[ (\bM^{(n)}_\ell)^\top \, \big| \, \cF^{(n)}_{\ell-1} \bigr]
    = \bzero$
 \ by \eqref{m}, thus we obtain \eqref{Mcond} and \eqref{Cov} in case
 \ $k \ne \ell$.

By \eqref{X} and \eqref{mean}, we conclude
 \[
   \bX^{(n)}_k - \EE(\bX^{(n)}_k)
   = \bm_{\bxi}^k (\bX^{(n)}_0 - \EE(\bX^{(n)}_0))
     + \sum_{j=1}^k \bm_{\bxi}^{k-j} \bM^{(n)}_j .
 \]
Now by \eqref{Cov},
 \begin{align*}
  \var(\bX^{(n)}_k)
  & = \bm_{\bxi}^k
      \EE\bigl[ (\bX^{(n)}_0 - \EE(\bX^{(n)}_0))
                (\bX^{(n)}_0 - \EE(\bX^{(n)}_0))^\top \bigr] (\bm_{\bxi}^\top)^k \\
  & \quad
      + \sum_{j=1}^k \sum_{\ell=1}^k
         \big(\bm_{\bxi}^\top \big) ^ {k-j}
        \EE \bigl[ \bM^n_j (\bM^n_\ell)^\top \bigr] \,
        (\bm_{\bxi}) ^ {k-\ell} \\
  & = \bm_{\bxi}^k \var(\bX^{(n)}_0) (\bm_{\bxi}^\top)^k
      + \sum_{j=1}^k
         \bm_{\bxi}^{k-j}
         \EE\bigl[ \bM^{(n)}_j (\bM^{(n)}_j)^\top \bigr] \,
         (\bm_{\bxi}^\top)^{k-j} .
 \end{align*}
Finally, using the expression in \eqref{Cov} for
 \ $\EE\bigl[ \bM^{(n)}_j (\bM^{(n)}_j)^\top \bigr]$ \ we obtain \eqref{var}.
\proofend

\begin{Lem}\label{EEX}
Under the assumptions of Theorem \ref{main} we have
 \begin{gather*}
  \EE(\|\bX^{(n)}_k\|) = \OO(k + n) , \qquad
  \EE(\|\bX^{(n)}_k\|^2) = \OO((k + n)^2) , \\
  \EE(\|\bM^{(n)}_k\|) = \OO((k + n)^{1/2}), \qquad
  \EE(\|\bM^{(n)}_k\|^4) = \OO((k + n)^2) .
 \end{gather*}
\end{Lem}

\noindent
\textit{Proof.}
By \eqref{mean},
 \[
   \|\EE(\bX^{(n)}_k)\|
   \leq \|\bm_{\bxi}^k\| \cdot \EE(\|\bX^{(n)}_0\|)
        + \sum_{j=0}^{k-1}
           \|\bm_{\bxi}^j\| \cdot \|\bm_{\bvare}\|
   \leq C_{\bm_{\bxi}} ( \sqrt{C} n + \|\bm_{\bvare}\| k ) ,
 \]
 where
 \begin{equation}\label{C}
  C_{\bm_{\bxi}}
  := \sup_{j \in \ZZ_+} \| \bm_{\bxi}^j \|
  < \infty, \qquad
  C
  := \sup_{n \in \NN} n^{-2} \EE(\| \bX^{(n)}_0 \|^2)
  < \infty ,
 \end{equation}
 since \eqref{rate} implies
 \ $C_{\bm_{\bxi}} \leq c_{\bm_{\bxi}} + \|\bPi_{\bm_\bxi}\|$. 
\ Hence, we obtain
 \ $\EE(\|\bX^{(n)}_k\|) \leq p \|\EE(\bX^{(n)}_k)\| = \OO(k + n)$. 

We have
 \begin{align*}
  \EE(\|\bM^{(n)}_k\|)
  & \leq \sqrt{\EE(\|\bM^{(n)}_k\|^2)}
    = \sqrt{\EE\bigl[ \tr (\bM^{(n)}_k (\bM^{(n)}_k)^\top) \bigr]}
    = \sqrt{\tr\bigl[ \bV_{\bvare} +  \EE(\bX_{k-1}^{(n)}) \odot \bV_\bxi \bigr]} \\
  & \leq \sqrt{\tr(\bV_{\bvare})}
         + \sqrt{\tr\bigl[ \EE(\bX_{k-1}^{(n)}) \odot \bV_\bxi \bigr]} ,
 \end{align*}
 hence we obtain \ $\EE(\|\bM^{(n)}_k\|) = \OO((k + n)^{1/2})$ \ from
 \ $\EE(\|\bX^{(n)}_k\|) = \OO(k + n)$. 

We have
 \[
   \EE(\|\bX^{(n)}_k\|^2) = \EE\bigl[\tr(\bX^{(n)}_k(\bX^{(n)}_k)^\top)\bigr]
   = \tr(\var(\bX^{(n)}_k))
     + \tr\bigl[ \EE(\bX^{(n)}_k) \EE(\bX^{(n)}_k)^\top \bigr] ,
 \]
 where
 \ $\tr\bigl[ \EE(\bX^{(n)}_k) \EE(\bX^{(n)}_k)^\top \bigr]
    = \|\EE(\bX^{(n)}_k)\|^2 \leq \bigl[\EE(\|\bX^{(n)}_k\|)\bigr]^2
    = \OO((k + n)^2)$.
\ Moreover, \ $\tr(\var(\bX^{(n)}_k)) = \OO((k + n)^2)$.
\ Indeed, by \eqref{var} and \eqref{C},
 \begin{align*} 
  \|\var(\bX^{(n)}_k)\|
  & \leq \sum_{j=0}^{k-1}
          \left( \|\bV_{\bvare}\|
                 + \|\bV_\bxi\| \cdot \|\bm_{\bxi}^{k-j-1}\|
                   \cdot \EE(\|\bX^{(n)}_0\|) \right)
          \| \bm_{\bxi}^j \|^2 \\
  & \quad
         + \|\var(X^n_0)\| \cdot \|\bm_{\bxi}^k\|^2
         + \|\bm_{\bvare}\| \cdot \|\bV_\bxi\|
           \sum_{j=0}^{k-2}
            \|\bm_{\bxi}^j\|^2
            \sum_{\ell=0}^{k-j-2}
             \|\bm_{\bxi}^\ell\| \\
  & \leq \left( \|\bV_{\bvare}\|
                + C_{\bm_{\bxi}} \|\bV_\bxi\| \cdot \EE(\|\bX^{(n)}_0\|) \right)
         C_{\bm_{\bxi}}^2 k \\
  & \quad
         + \big( \EE(\|\bX^{(n)}_0\|^2) + \bigl[\EE(\|X^n_0\|)\bigr]^2 \big)
           C_{\bm_{\bxi}}^2
         + C_{\bm_{\bxi}}^3 \|\bm_{\bvare}\| \cdot \|\bV_\bxi\| k^2 ,
 \end{align*}
 where \ $\|\bV_\bxi\| := \sum_{i=1}^p \|\bV_{\bxi_i}\|$, \ hence we obtain
 \ $\EE(\|\bX^{(n)}_k\|^2) = O((k + n)^2)$.

By \eqref{Mdeco}, 
 \[
   \|\bM^{(n)}_k\|
   \leq \|\bvare^{(n)}_k - \EE(\bvare^{(n)}_k)\|
        + \sum_{i=1}^p
           \Biggl\| \sum_{j=1}^{X_{k-1,i}^{(n)}}
                    ( \bxi^{(n)}_{k,j,i} - \EE(\bxi^{(n)}_{k,j,i}) ) \Biggr\| ,
 \]
 hence
 \[
   \EE(\|\bM^{(n)}_k\|^4)
   \leq (p+1)^3 \EE(\|\bvare^{(1)}_1 - \EE(\bvare^{(1)}_1)\|^4)
        + (p+1)^3
          \sum_{i=1}^p
           \EE\Biggl(\Biggl\| \sum_{j=1}^{X_{k-1,i}^{(n)}}
                             ( \bxi^{(n)}_{k,j,i} - \EE(\bxi^{(n)}_{k,j,i}) )
                     \Biggr\|^4
              \Biggr) .
 \]
Here
 \begin{align*}
  \EE\Biggl(\Biggl\| \sum_{j=1}^{X_{k-1,i}^{(n)}}
                      ( \bxi^{(n)}_{k,j,i} - \EE(\bxi^{(n)}_{k,j,i}) ) \Biggr\|^4
            \Biggr)
  &= \EE\Biggl[\Biggl(\sum_{\ell=1}^p
                       \Biggl( \sum_{j=1}^{X_{k-1,i}^{(n)}}
                               ( \xi^{(n)}_{k,j,i,\ell} - \EE(\xi^{(n)}_{k,j,i,\ell}) )
                       \Biggr)^2
              \Biggr)^2\Biggr] \\
  &\leq p\sum_{\ell=1}^p
          \EE\Biggl[\Biggl(\sum_{j=1}^{X_{k-1,i}^{(n)}}
                           ( \xi^{(n)}_{k,j,i,\ell} - \EE(\xi^{(n)}_{k,j,i,\ell}) )
                    \Biggr)^4
             \Biggr] ,
 \end{align*}
 where
 \begin{multline*}
  \EE\Biggl[\Biggl(\sum_{j=1}^{X_{k-1,i}^{(n)}}
                  ( \xi_{k,j,i,\ell} - \EE(\xi_{k,j,i,\ell}) ) \Biggr)^4
           \, \Bigg| \, \cF_{k-1}^{(n)}\Biggr] \\
  = X_{k-1,i}^{(n)} \EE[( \xi^{(1)}_{1,1,i,\ell} - \EE(\xi^{(1)}_{1,1,i,\ell}) )^4]
     + X_{k-1,i}^{(n)} (X_{k-1,i}^{(n)} - 1)
       \bigl(\EE[( \xi^{(1)}_{1,1,i,\ell} - \EE(\xi^{(1)}_{1,1,i,\ell}) )^2]\bigr)^2 
 \end{multline*}
 with
 \ $\bigl(\EE[( \xi^{(1)}_{1,1,i,\ell} - \EE(\xi^{(1)}_{1,1,i,\ell}) )^2]\bigr)^2
    \leq \EE[( \xi^{(1)}_{1,1,i,\ell} - \EE(\xi^{(1)}_{1,1,i,\ell}) )^4]$,
 \ hence
 \[
   \EE\Biggl[\Biggl(\sum_{j=1}^{X_{k-1,i}^{(n)}}
                     ( \xi^{(n)}_{k,j,i,\ell} - \EE(\xi^{(n)}_{k,j,i,\ell}) ) \Biggr)^4
      \Biggr]
   \leq \EE[( \xi^{(1)}_{1,1,i,\ell} - \EE(\xi^{(1)}_{1,1,i,\ell}) )^4]
        \EE[(X^{(n)}_{k-1,i})^2]
 \]
Consequently, \ $\EE(\|\bX^{(n)}_k\|^2) = \OO((k + n)^2)$ \ implies
 \ $\EE(\|\bM^{(n)}_k\|^4) = \OO((k + n)^2)$.
\proofend

Next we recall a result about convergence of random step processes towards a
 diffusion process, see Isp\'any and Pap \cite[Corollary 2.2]{IspPap}.

\begin{Thm}\label{Conv2DiffThm}
Let \ $\bgamma : \RR_+ \times \RR^p \to \RR^{p \times r}$ \ be a continuous
 function.
Assume that uniqueness in the sense of probability law holds for the SDE
 \begin{equation}\label{SDE}
  \dd \, \bcU_t
  = \gamma (t, \bcU_t) \, \dd \bcW_t ,
  \qquad t \in \RR_+,
 \end{equation}
 with initial value \ $\bcU_0 = \bu_0$ \ for all \ $\bu_0 \in \RR^p$, \ where
 \ $(\bcW_t)_{t \in \RR_+}$ \ is an $r$-dimensional standard Wiener process.
Let \ $\bmu$ \ be a probability measure on \ $(\RR^p, \cB((\RR^p))$, \ and let
 \ $(\bcU_t)_{t \in \RR_+}$ \ be a solution of \eqref{SDE} with initial
 distribution \ $\bmu$.

For each \ $n \in \NN$, \ let \ $(\bU^{(n)}_k)_{k \in \ZZ_+}$ \ be a sequence of
 $p$-dimensional martingale differences with respect to a filtration
 \ $(\cF^{(n)}_k)_{k \in \ZZ_+}$.
\ Let
 \[
   \bcU^{(n)}_t := \sum_{k=0}^{\nt} \bU^{(n)}_k \, ,
   \qquad t \in \RR_+, \quad n \in \NN .
 \]
Suppose \ $\EE \big( \|\bU^{(n)}_k\|^2 \big) < \infty$ \ for all
 \ $n, k \in \NN$, \ and \ $\bU^{(n)}_0 \distr \bmu$.
\ Suppose that, for each \ $T > 0$,
 \begin{enumerate}
  \item[\textup{(i)}]
   $\sup\limits_{t\in[0,T]}
     \left\| \sum\limits_{k=1}^{\nt}
              \EE\Bigl[ \bU^{(n)}_k (\bU^{(n)}_k)^\top \mid \cF^{(n)}_{k-1} \Bigr]
             - \int_0^t
                \bgamma(s,\bcU^{(n)}_s) \bgamma(s,\bcU^{(n)}_s)^\top
                \dd s \right\|
         \stoch 0$,\\
  \item[\textup{(ii)}]
   $\sum\limits_{k=1}^{\lfloor nT \rfloor}
     \EE \big( \|\bU^{(n)}_k\|^2 \bbone_{\{\|\bU^{(n)}_k\| > \theta\}}
                 \bmid \cF^{(n)}_{k-1} \big)
    \stoch 0$
   \ for all \ $\theta > 0$,
 \end{enumerate}
 where \ $\stoch$ \ denotes convergence in probability.
Then \ $\bcU^{(n)} \distr \bcU$ \ as \ $n \to \infty$.
\end{Thm}

Now we recall a version of the continuous mapping theorem.

For functions \ $f$ \ and \ $f_n$, \ $n \in \NN$, \ in \ $\DD(\RR_+, \RR^p)$,
 \ we write \ $f_n \lu f$ \ if \ $(f_n)_{n \in \NN}$ \ converges to \ $f$
 \ locally uniformly, i.e., if \ $\sup_{t \in [0,T]} \|f_n(t) - f(t)\| \to 0$ \ as
 \ $n \to \infty$ \ for all \ $T > 0$.
\ For measurable mappings \ $\Phi : \DD(\RR_+, \RR^p) \to \DD(\RR_+, \RR^q)$
 \ and \ $\Phi_n : \DD(\RR_+, \RR^p) \to \DD(\RR_+, \RR^q)$, \ $n \in \NN$, \ we
 will denote by \ $C_{\Phi, (\Phi_n)_{n \in \NN}}$ \ the set of all functions
 \ $f \in \CC(\RR_+, \RR^p)$ \ for which \ $\Phi_n(f_n) \to \Phi(f)$ \ whenever
 \ $f_n \lu f$ \ with \ $f_n \in \DD(\RR_+, \RR^p)$, \ $n \in \NN$.

\begin{Lem}\label{Conv2Funct}
Let \ $(\bcU_t)_{t \in \RR_+}$ \ and \ $(\bcU^{(n)}_t)_{t \in \RR_+}$, \ $n \in \NN$,
 \ be \ $\RR^p$-valued stochastic processes with c\`adl\`ag paths such that
 \ $\bcU^{(n)} \distr \bcU$.
\ Let \ $\Phi : \DD(\RR_+, \RR^p) \to \DD(\RR_+, \RR^q)$ \ and
 \ $\Phi_n : \DD(\RR_+, \RR^p) \to \DD(\RR_+, \RR^q)$, \ $n \in \NN$, \ be
 measurable mappings such that there exists \ $C \subset C_{\Phi,(\Phi_n)_{n\in\NN}}$
 \ with \ $C \in \cD_\infty(\RR_+, \RR^p)$ \ and \ $\PP(\bcU \in C) = 1$.
\ Then \ $\Phi_n(\bcU^{(n)}) \distr \Phi(\bcU)$.
\end{Lem}

Lemma \ref{Conv2Funct} can be considered as a consequence of Theorem 3.27 in
 Kallenberg \cite{K}, and we note that a proof of this lemma can also be found
 in Isp\'any and Pap \cite[Lemma 3.1]{IspPap}.

\end{document}